\documentclass{article}
\usepackage{amsthm}

\newtheorem{Lemma}{Lemma}

\newtheorem{Property}{Property}[section]

\newcommand{\bincoeff}[2]{\mbox{$\left(\!\begin{array}{c}{#1}\\{#2}\end{array}\!\right)$}}

\newcommand{\legendre}[2]{\mbox{$\left(\!\!\begin{array}{c}{\underline{{~\raisebox{0.25em}{\mbox{$#1$}}~}}}\\{\mbox{$#2$}}\end{array}\!\!\right)$}}

\newcommand{\notdiv}{\mbox{$\hspace*{4pt}|\hspace*{-3.75pt}/~$}}
\newcommand{\mod}[1]{\mbox{$\!\!\pmod{#1}$}}
\newcommand{\myProof}{\proof}                                  

\begin{document}


\title{On Hunting for Taxicab Numbers}

\author{P. Emelyanov \\
        Institute of Informatics Systems \\
        6 avenue Lavrentiev, 630090, Novosibirsk, Russia \\
        e-mail: emelianov@iis.nsk.su}

 \maketitle
\begin{abstract}

In this article, we make use of some known method to investigate
some properties of the numbers represented as sums of two equal
odd powers, i.e., the equation $x^n+y^n=N$\/ for $n\ge3$. It was
originated in developing algorithms to search new taxicab numbers
(i.e., naturals that can be represented as a sum of positive cubes
in many different ways) and to verify their minimality. We discuss
properties of diophantine equations that can be used for our
investigations. This techniques is applied to develop an algorithm
allowing us to compute new taxicab numbers (i.e., numbers
represented as sums of two positive cubes in $k$ different ways),
for $k=7\ldots14$\/.

\end{abstract}

\section*{Introduction}

This work was originated in searching new so--called {\em taxicab
numbers}, i.e., naturals $T_k$\/ that can be
represented/decomposed as/into a sum of positive cubes in $k$
different ways, and verifying their minimality. We made use of
some known method to investigate properties of the cubic equation
that could help us to find new taxicab numbers.

Already Fermat proved that numbers expressible as a sum of two
cubes in $n$\/ different ways exist for any $n$. But still finding
taxicab numbers and proving their minimality are hard
computational problems. Whereas the first nontrivial taxicab
number $T_2=1729$\/ became widely--known in 1917 thanks to
Ramanujan and Hardy, next ones were only found with help of
computers: $T_3=87539319$ (J. Leech, 1957), $T_4=6963472309248$
(E. Rosenstiel, J.A. Dardis, and C.R. Rosenstiel, 1991),
$\mbox{\bf\em W}_5=T_5=48988659276962496$ (D. Wilson, 1997,
\cite{Wilson-jis-1999}). It is known that these numbers are
minimal. For $\mbox{\bf\em R}_6=T_6=24153319581254312065344$ (R.L.
Rathbun, 2002) as well as for next discovered taxicab numbers it
is unknown.

In January--September 2006 the author computed
$T_7=139^3\mbox{\bf\em R}_6$, $T_8=727^3T_7$, $T_9=4327^3T_8$,
$T_{10}=38623^3T_9$, and $T_{11}=45294^3T_{10}$. At the end of
2006 the author learned about the results of C. Boyer
\cite{Boyer-2006} who established smaller $T_7,\ldots,T_{11}$\/
and first $T_{12}$\/ in December 2006. At the begin of 2007 the
author computed $T_{13}$\/ and $T_{14}$.

The article is organized as follows. We start with putting the
equation in a new form. Next, we deduce simple properties of the
equation of interest based on this presentation. At the end, we
present a new algorithm to compute taxicab numbers which we used
to find new ones.

\section{Common Properties}

We are interested in the problem of representations (also called
decompositions) of numbers as the sums of two positive odd
$n$-powers; i.e., solvability of the equation
\begin{equation}\label{OriginalEquation}
x^n+y^n=N
\end{equation}
in positive integers. A solution of this equation is also called a
representation or a decomposition of the number $N$. The equation
of interest is too ``smooth'' in its original form. We want to
make it ``uneven''. We are going to consider this equation in the
following $m\pm{}h$-form ($m\neq h>0$)
\begin{equation}\label{TheEquation}
(m-h)^n+(m+h)^n=N
\end{equation}
which is not an infrequent guest in number--theoretical proofs.

Although only even numbers can be directly represented in this
way, there is a simple transformation that allows us to treat this
equation for odd $N$\/ as well. In fact, any pair $(x,y)$\/
consisting of even and odd integers can be represented as
$(t-s-1,t+s)$. If $N$\/ is odd, we write
\[
(t-s-1)^n+(t+s)^n=N.
\]
Multiplying both sides by $2^n$\/ we can put the previous equation
into the form
\begin{equation}\label{TheEqOdd}
((2t-1)-(2s+1))^n+((2t-1)+(2s+1))^n=2^nN
\end{equation}
and, then some extra steps are needed to obtain representations of
$N$\/ itself. For the exponent 3, the least odd number $N$\/ for
which $2^3N$\/ yields a not only proper two cubes representation
is 513:
\[
2^3 513=2^3\left(1^3+8^3\right)=(12-3)^3+(12+3)^3=4104.
\]
Notice that 4104 is the least even number represented as a sum of
two cubes in two different ways. Next, assume $N$\/ to be even if
we do not explicitly state the contrary.

We are interested in any prime powers, although sometimes it is
sufficient that they are odd only. Such representations for odd
powers are closely related to divisors of the numbers of interest.



We shall refer to $m$\/ as a {\em median}\/ of the corresponding
power representation and to $N_d$\/ as an integer quotient $N/d$\/
if it exists. We shall make use the following property (a simple
corollary of Quadratic Reciprocity Law) of odd prime divisors of
binary forms:
\begin{Property}\label{BinaryForm3}
\[
p ~|~ ax^2+by^2 ~\wedge~ \gcd(ax,by)=1 ~~\Longrightarrow~~
\legendre{ab}{p}=(-1)^\frac{p-1}{2}.
\]
In particular, for the binary form $u^2+3v^2$\/ the forbidden
divisors are
\[
\legendre{3}{p}\not=(-1)^\frac{p-1}{2}; \mbox{~~i.e.,~~}
p\equiv5,11\mod{12}.
\]
\end{Property}

Given $N$\/ and its divisors, by solving an $n-1$-order polynomial
equation     
\begin{equation}\label{ExpansionEQ}
(m-h)^n+(m+h)^n=
2m\left(%
\sum_{k=0}^{\frac{n-1}{2}}\bincoeff{n}{2k}m^{n-2k-1}h^{2k}%
\right)=N
\end{equation}
with respect to $h$, we can either ``easily'' find some
representation(s) of this number or prove that it is impossible.
Notice that in this polynomial $m$\/ and $h$\/ occur only in odd
and even powers, respectively.

We start the investigation by establishing the following simple
properties of {\bf Equation (\ref{TheEquation})}.

\begin{Lemma}\label{DivNModulo}
If $m$\/ is a median of some representation of $N$, then
\[
m\equiv N_2\mod{n}.
\]
If $n | N$, then also $n^2 | N$. If $n \notdiv{} N$, then
$N=2m(nt+1)$.
\end{Lemma}

\myProof
First, rewriting {\bf Equation (\ref{TheEquation})} in the form
\[
m^n+n\sum_{k=1}^{\frac{n-1}{2}}\frac{1}{n}\bincoeff{n}{2k}m^{n-2k}h^{2k}=N_2
\]
we can derive the modular equation $m^n\equiv N_2\mod{n}$. Next:
\begin{itemize}

\item By applying Fermat's Little Theorem we have the first
statement.

\item Because $n | N$, therefore also $n | m$, and this yields the
second statement.

\item Because $n \notdiv{} N$, then also $n \notdiv{} m$. By
applying Fermat's Little Theorem to $m^{n-1}\equiv
N_{2m}\mod{n}$\/ we have the third statement.

\end{itemize}
\qed

Because $h$\/ is ranged in $(0,m)$\/ it is easy to establish
\begin{Lemma}\label{DivBounds}
A necessary condition for $N$\/ to have a representation as the
sum of $n$-powers is
\[
 \exists m|N ~~:~~ \sqrt[n]{\frac{N}{2^n}}~<~m~<~\sqrt[n]{\frac{N}{2}}.
\]
\end{Lemma}
\noindent Obviously, the number of such representations does not
exceed the number of divisors of $N$\/ satisfying this condition
(see also {\bf Lemma \ref{TaxicabLowerBound}}).

{\bf Lemmas 1} and {\bf2} allow us to estimate numbers being the
sum of two odd powers higher than 2 in $k$\/ ways. If a number has
two different representations for the power $n$, then the medians
$m_1, m_2$\/ corresponding to them also satisfy the congruence
$m_1\equiv m_2\mod{n}$. Because
$\sqrt[n]{N/2^n}+n(k-1)\leq\sqrt[n]{N/2}$, we have the following
properties of generalized taxicab numbers
\begin{Lemma}\label{TaxicabLowerBound}
If number $T(n,k)$, $k>1$, represented as the sum of two
$n$-powers in $k$\/ ways is even, then it has at least $k$\/
divisors in the range $(\sqrt[n]{N/2^n},\sqrt[n]{N/2})$\/ and the
following lower bound holds
\[
T(n,k)\geq2\left(\frac{2n}{2-\sqrt[n]{2}}\right)^n(k-1)^n
\]

\end{Lemma}

\noindent This bound is far from optimal due to a quite
conservative assumption about the gaps between medians. This is a
subject of further investigation. Recall that only wide-known
theoretical bound for $T(3,k)=T_k$\/ is Silverman's result
\cite{Silverman-jlms-1983} that  describes its logarithmic
behavior:
\[
\log T_k=o(k^{r+2/r}),
\]
where $r$\/ is the highest rank of {\bf Equation
(\ref{OriginalEquation})}. The highest rank known now is 5.

When there are "too many" taxicab medians they cannot be relative
prime because all of them are divisors. Hence they share common
divisors. In particular, for taxicab medians $m_1<\ldots<m_k$\/
the following inequality holds:
\[
\mbox{lcm}(m_1,\ldots,m_k)\leq (2m_1)^n.
\]

The cubic equation in the form $m^2+3h^2=N_{2m}$\/ provides a way
to derive parameterizations of the two cubes representation
problem\footnote{Here we treat independently the median $m$\/ and
its co-factor $N_{2m}$; therefore this does not cover general
cases.}. We mention only those of them that relate to the taxicab
numbers problem. It arises when $N_{2m}$\/ is a cube and this case
is connected to the well-known problem of the decomposition of
numbers into two rational cubes (positive or not) which was
investigated by Fermat, Euler, Sylvester, and other researchers.

G\'erardin proved \cite[Chapter XX]{Dickson-1999} that all
solutions of $u^2+3v^2=w^3$\/ with $\gcd(u,v)=1$\/ are generated
by \label{Gerardin}
\[
(t^3-9ts^2)^2+3(3t^2s-3s^3)^2=(t^2+3s^2)^3.
\]
We have
\[
(t^3-3t^2s-9ts^2+3s^3)^3+(t^3+3t^2s-9ts^2-3s^3)^3=2(t^3-9ts^2)(t^2+3s^2)^3,
\]
and next
\[
\left(\frac{t^3-3t^2s-9ts^2+3s^3}{t^2+3s^2}\right)^3+
\left(\frac{t^3+3t^2s-9ts^2-3s^3}{t^2+3s^2}\right)^3=2t(t-3s)(t+3s).
\]
So, if the diophantine equation
$2t^3-18ts^2\pm{}Nr^3=0$\/ is solvable%
\footnote{Euler's solution of the two rational cubes problem is
slightly different.}, then $N$\/ is decomposable.

This can be simplified into one-parametric examples as follows
\[
\left(\frac{w^3+3w^2-6w+1}{3(w^2-w+1)}\right)^3-
\left(\frac{w^3-6w^2+3w+1}{3(w^2-w+1)}\right)^3=w(w-1),
\]
and
\[
\pm{}\left(\frac{8w^9\pm{}24w^6+6w^3\mp{}1}{3w(4w^6\pm{}2w^3+1)}\right)^3\mp{}
\left(\frac{8w^9\mp{}12w^6-12w^3\mp{}1}{3w(4w^6\pm{}2w^3+1)}\right)^3=4w^3\pm{}2.
\]
Also, the substitution $t-3s=u^2v, t+3s=uv^2$\/  gives
\[
\left(\frac{u^3+6u^2v+3uv^2-v^3}{3(u^2+uv+v^2)}\right)^3+
\left(\frac{v^3+6v^2u+3vu^2-u^3}{3(u^2+uv+v^2)}\right)^3=uv(u+v)
\]
which provides the following parametrization of the sum of two
integer powers
\[
\left(\frac{p^9+6p^6q^3+3p^3q^6-q^9}{3pq(p^6+p^3q^3+q^6)}\right)^3+
\left(\frac{q^9+6q^6p^3+3q^3p^6-p^9}{3pq(p^6+p^3q^3+q^6)}\right)^3=p^3+q^3.
\]

Catalan's parametrization
\[
\left({\mbox{\small$\frac12$}}\,(t+s)(t-2s)(s-2t)\right)^2+
3\left({\mbox{\small$\frac32$}}\,ts(t-s)\right)^2=
\left(t^2-ts+s^2\right)^3
\]
leads us to another rational cubes identity
\[
\left(\frac{t^3-3t^2s+s^3}{t^2-ts+s^2}\right)^3+
\left(\frac{t^3-3ts^2+s^3}{t^2-ts+s^2}\right)^3=(t+s)(2s-t)(s-2t).
\]
The substitution $2s-t=u^2v, s-2t=uv^2$\/ gives the following
identity
\[
\left(\frac{u^3+3u^2v-6uv^2+v^3}{3(u^2-uv+v^2)}\right)^3-
\left(\frac{u^3-6u^2v+3uv^2+v^3}{3(u^2-uv+v^2)}\right)^3=uv(u-v)
\]
which provides the following parametrization of the sum of two
integer powers
\[
\left(\frac{p^9+3p^6q^3-6p^3q^6+q^9}{3pq(p^6-p^3q^3+q^6)}\right)^3-
\left(\frac{p^9-6p^6q^3+3p^3q^6+q^9}{3pq(p^6-p^3q^3+q^6)}\right)^3=p^3-q^3.
\]

It is easy to note that these parameterizations of the sum and the
difference of two integer cubes also give parameterizations to the
diophantine equation $X^3+Y^3=S^3+T^3$. Euler's parametric
solution to $X^3+Y^3=S^3+T^3$\/ is
\[
\begin{array}{lll}
X = w (1-(u -3v)(u^2+ 3v^2))   &~~~& Y = w ((u + 3v)(u^2+ 3v^2)-1) \\
S = w ((u + 3v)-(u^2+ 3v^2)^2) &   & T = w ((u^2+ 3v^2)^2+ (3v-u))
\end{array}
\]


Finally, we mention some properties of the equation of interest
that can be used to investigate taxicab numbers. Sometimes we can
improve the congruence of {\bf Lemma \ref{DivNModulo}}:

\begin{Lemma}
If $(m-h)^p+(m+h)^p=N$,~ $\gcd(m,h)=1$\/,~ $m\not\equiv
h\mod{2}$,~ then

\[
p=3  ~~\Rightarrow~~ m\equiv N_2\mod{12}
\]
\[
p=5  ~~\Rightarrow~~ m\equiv N_2\mod{20}.
\]

\end{Lemma}

\myProof We write down $2m^3+6mh^2=N$\/ as $m^2-h^2+4h^2=N_{2m}$\/
and $2m^5+20m^3h^2+10mh^4=N$\/ as $5m(m^2+h^2)^2-4m^5=N_2$.
Considering these equations by modulo 4 we conclude $m\equiv
N_2\mod{4}$. Combining this congruence with the congruence from
{\bf Lemma \ref{DivNModulo}} we obtain these lemma statements.

\qed

The forbidden divisors condition for two-squares representation is
well known since Fermat's work. For cubic and quintic equations
there are analogies which follow from {\bf Property
\ref{BinaryForm3}}:

\begin{Lemma}\label{ForbiddenDivisors}
Necessary conditions for $N$\/ to have a cubic/quintic
representation with $\gcd(m,h)=1$\/ are the following:
\begin{enumerate}
\item It has no prime divisors of forms $12t+5$\/ and $12t+11$\/
(the cubic case) or of forms $10t\pm{}1$\/ (the quintic case), or

\item If such divisors exist, then all of them are factors of the
median.
\end{enumerate}
\end{Lemma}

\noindent{\bf Remark.}~ In view of the cubic case of {\bf Lemma
\ref{ForbiddenDivisors}}, we can mention the results of Euler et
al for the divisors of numbers in the form $u^2+3v^2$: all prime
divisors have the same form $\alpha^2+3\beta^2$.
\vskip\baselineskip

\section{New Taxicab/Cabtaxi Numbers}

Before we discuss cubic taxicab numbers, we briefly consider the
equation $x^5+y^5=u^5+v^5$. No such number is known within the
range up to $1.05\cdot10^{26}$. We have not yet found any, but we
found some solution in Gaussian integers:
\[
\begin{array}{c}
\left(t^2+s^2-(t^2-2ts-s^2)\imath\right)^5+\left(t^2+s^2+(t^2-2ts-s^2\right)\imath)^5=\\
\left(t^2+s^2-(t^2+2ts-s^2)\imath\right)^5+\left(t^2+s^2+(t^2+2ts-s^2\right)\imath)^5=\\
-8(t^2+s^2)(t^4-2t^3s-6t^2s^2+2ts^3+s^4)(t^4+2t^3s-6t^2s^2-2ts^3+s^4)
\end{array}
\]
The least such positive number is
$3800=(5-\imath)^5+(5+\imath)^5=(5-7\imath)^5+(5+7\imath)^5$.

The observation that $T_6=79^3T_5$\/ stirs up our interest in
searching for new taxicab numbers $T_k$\/ in the same way. The
usual definition of taxicab numbers is equipped with a condition
that they are minimal. But for brevity we designate all multi-ways
representable numbers as taxicab numbers.
Even an open question%
\footnote{%
C. Calude et al \cite{CaludeCaludeDinneen-jucs-2003} (with an
update \cite{CaludeCaludeDinneen-CDMTCS-2005}) stated that the
minimality of $T_6$\/ can be confirmed with the probability
$>0.99$\/ but G. Martin criticized their considerations in
Mathematical Reviews MR2149410 (2006a:11175).%
}%
~about the minimality of $T_6$\/ does not matter. To compute some
$k+1$--way representable number we can try  any $k$-way
representable number. Our approach can produce non-minimal
numbers, but such numbers can be used to check their minimality or
to search for smaller ones. We believe that this median--based
approach reducing the length of tested numbers in three times
allows us to check the minimality of $T_6$\/ and $T_7$.

Notice that Wilson \cite{Wilson-jis-1999} used similar ideas
(cubic multipliers) to find 5--way representable number
$\mbox{\bf\em W}_5=48988659276962496$\/ in 1997 but his approach
is more expensive even for small numbers. During this search a
six-way example was also detected. Inspired by Wilson's approach
in 2002 R. L. Rathbun \cite{Rathbun-2002} presented the smaller
candidate
\[
\mbox{\bf\em R}_6=79^3\,\mbox{\bf\em W}_5=24153319581254312065344.
\]
Rathbun also mentioned multipliers $139$\/ and $727$\/ giving
other examples of six-way representable numbers. Our approach
demonstrates that they appear in multipliers of $T_9$\/ and
$T_{11}$, respectively.

In the first version of this article (December 2006) we described
a modification of our algorithm that  produces some taxicab
numbers. In January--September 2006 with help of this algorithm we
computed $T_7=139^3\mbox{\bf\em R}_6$, $T_8=727^3T_7$,
$T_9=4327^3T_8$, $T_{10}=38623^3T_9$, and $T_{11}=45294^3T_{10}$.
At that moment we learned about results of C. Boyer
\cite{Boyer-2006} who established smaller $T_7,\ldots,T_{11}$\/
and first $T_{12}$\/ in December 2006. Unfortunately he has not
yet published details of his algorithm. Our renewed algorithm,
given later in this article, produces the same numbers. Also, for
the first time we found $T_{13}$\/ and $T_{14}$.

The main idea of our approach is not too surprising. If we know
some $k$--way representable number $T_k$, then we can try to find
$T_{k+1}$\/ in the form $\mu^3\,T_k$. If $m_1,\ldots,m_k$\/ are
medians of the representations of $T_k$, then medians of the
representations of $T_{k+1}$\/ are
$\mu{}m_1,\ldots,\mu{}m_k,d^\prime{}d$\/ where
$d^\prime\in\mbox{\rm divisors}(\mu^3)$\/ (the first version of
the algorithm uses only $d^\prime=1$) and $d\in\mbox{\rm
divisors}(T_k)$. A simple observation is that the multiplier of
interest does not exceed $2T_k^{2/3}$.

The iterative procedure formalizing this idea and using the
properties of the equation is the following:

\begin{description}
\item[$\bullet$~~~] Create an ordered array $D$\/ of all divisors
of $T_k$\/ excluding known too small divisors, i.e., less than
$\sqrt[3]{T_k/4}$.

\item[$\bullet$~~~] For multipliers $M$\/ from 2 to $\lfloor
2T_k^{2/3}\rfloor$\/ do

\begin{description}

 \item[$\bullet$~~~] Let $N=M^3T_k$;

 \item[$\bullet$~~~] For $\mu\in\mbox{divisors}(M^3)$\/ do

\begin{description}

\item[$\bullet$~~] Using dichotomic search, find a range of $D$\/
where the divisors satisfying {\bf Lemma \ref{DivBounds}} for
$\frac1\mu{}N$\/ are located;

\item[$\bullet$~~] Within this range for divisors $d$\/ such that
$\mu{}d\equiv \frac12N\!\mod{3}$\/ do: if the value
$(\frac{1}{2\mu{}d}N-(\mu{}d)^2)/3$\/ is a perfect square, then
$\mu{}d$\/ is the $k+1$\/ median and therefore $N$\/ is $T_{k+1}$.
Otherwise continue.

\end{description}
\end{description}
\end{description}

A set of all divisors of $T_k$\/ may be space-consuming. To avoid
the explicit computation of this set we used the following trick.
A taxicab number $T_k$\/ is a product $M\!\cdot T_s$\/ where
$T_s$\/ is a ``seed'', i.e., a small taxicab number with an easily
computed set of divisors and $M=(\mu_{s+1}\cdots\mu_k)^3$.
Evidently $M=1$\/ for $T_{k+1}=T_{s+1}$. Thus computing
$T_{k+1}$\/ we split the loop iterating through all divisors of
$T_k$\/ into two nested loops: the outer loop iterating through
all divisors of $M$\/ and the inner one iterating through those
divisors of $T_s$\/ such that product of the first iterator, the
second iterator, and some divisor of the current cubic multiplier
satisfies {\bf Lemma \ref{DivBounds}}.

Choice of the seed $T_s$\/ affects the space used by the
algorithm. We used $\mbox{\bf\em W}_5$\/ to compute new $T_k$\/
for $k=7\ldots12$. But for the next numbers, cardinality of the
divisor set for $M$\/ exceeds one for $T_s$\/ more and more. To
balance the cardinalities of these sets we took greater seeds.

{\begin{center}
\begin{tabular}{||r|r|r|r||}
\hline\hline
Ways & Seed & Multiplier & Time~~~~~  \\
\hline\hline

7  & 5 & 101 & 58 s.   \\
8  & 5 & 127 & 5 m. 1 s.   \\
9  & 5 & 139 & 18 m. 47 s. \\
10 & 5 & 377 & 4 h. 8 m. \\
11 & 5 & 727  & 123 h. 20 m. \\
12 & 5 & 2971 & 152 d. \\
13 & 6 & 4327 & 21 h. 8 m.$^{a)}$ \\
14 & 6 & 7549 & 23 m. 39 s.$^{b)}$ \\

\hline\hline
\end{tabular}

{\small

\vskip\baselineskip\flushleft{$^{a)}$ To compute this number we
examined only prime multipliers great than 2971.}

\vspace*{-1em}\flushleft{$^{b)}$ To compute this number we
examined only this multiplier.}

}

\vskip\baselineskip {\bf Table 1.} Computational results.

\end{center}}

{\bf Table 1.} represents multipliers producing new taxicab
numbers. In {\bf APPENDIX A} we give these numbers themselves and
their decompositions.

Also, we found that all of our taxicab numbers $T(3,k)$\/ are {\em
cabtaxi} (i.e., without the restriction on the cubes of the
decomposition to be positive) numbers $C(3,k+2)$. Surprisingly the
multiplier 5 gives cabtaxi numbers of higher orders:
$5^3T(3,k)=C(3,k+4)$. We checked this property for $k=6\ldots12$.


\section*{Final Remark}

In September 2007 we learned about new results of C. Boyer who
established new taxicab numbers for $n=13\ldots19$\/ and cabtaxi
numbers for $n=10\ldots30$. Boyer's article is going to be
published in a mathematical magazine.


\providecommand{\bysame}{\leavevmode\hbox to3em{\hrulefill}\thinspace}
\providecommand{\MR}{\relax\ifhmode\unskip\space\fi MR }
\providecommand{\MRhref}[2]{%
  \href{http://www.ams.org/mathscinet-getitem?mr=#1}{#2}
}
\providecommand{\href}[2]{#2}

\section*{APPENDIX A. Taxicab numbers decompositions}

\noindent$T_7=101^3\,\mbox{\bf\em
R}_6=24885189317885898975235988544$:
\[
\begin{array}{rcrc}
58798362^3    & + &  2919526806^3 & =  \\
309481473^3   & + &  2918375103^3 & =  \\
459531128^3   & + &  2915734948^3 & =  \\
860447381^3   & + &  2894406187^3 & =  \\
1638024868^3  & + &  2736414008^3 & =  \\
1766742096^3  & + &  2685635652^3 & =  \\
1847282122^3  & + &  2648660966^3 &   \\
\end{array}
\]

\noindent$T_8=127^3\,T_7$=50974398750539071400590819921724352:
\[
\begin{array}{rcrc}
7467391974^3    & + & 370779904362^3 & =  \\
39304147071^3   & + & 370633638081^3 & =  \\
58360453256^3   & + & 370298338396^3 & =  \\
109276817387^3  & + & 367589585749^3 & =  \\
208029158236^3  & + & 347524579016^3 & =  \\
224376246192^3  & + & 341075727804^3 & =  \\
234604829494^3  & + & 336379942682^3 & =  \\
288873662876^3  & + & 299512063576^3 &
\end{array}
\]

\noindent$T_9=139^3\,T_8=136897813798023990395783317207361432493888$:
\[
\begin{array}{rcrc}
1037967484386^3  & + & 51538406706318^3 & = \\
4076877805588^3  & + & 51530042142656^3 & = \\
5463276442869^3  & + & 51518075693259^3 & = \\
8112103002584^3  & + & 51471469037044^3 & = \\
15189477616793^3 & + & 51094952419111^3 & = \\
28916052994804^3 & + & 48305916483224^3 & = \\
31188298220688^3 & + & 47409526164756^3 & = \\
32610071299666^3 & + & 46756812032798^3 & = \\
40153439139764^3 & + & 41632176837064^3 &  \\
\end{array}
\]

\noindent$T_{10}=377^3\,T_9=7335345315241855602572782233444632535674275447104$:
\[
\begin{array}{rcrc}
391313741613522^3   & + & 19429979328281886^3 & = \\
904069333568884^3   & + & 19429379778270560^3 & = \\
1536982932706676^3  & + & 19426825887781312^3 & = \\
2059655218961613^3  & + & 19422314536358643^3 & = \\
3058262831974168^3  & + & 19404743826965588^3 & = \\
5726433061530961^3  & + & 19262797062004847^3 & = \\
10901351979041108^3 & + & 18211330514175448^3 & = \\
11757988429199376^3 & + & 17873391364113012^3 & = \\
12293996879974082^3 & + & 17627318136364846^3 & = \\
15137846555691028^3 & + & 15695330667573128^3 &  \\
\end{array}
\]

\noindent$T_{11}=727^3\,T_{10}=2818537360434849382734382145310807703728251895897826621632$:
\[
\begin{array}{rcrc}
284485090153030494^3  & + & 14125594971660931122^3 & = \\
657258405504578668^3  & + & 14125159098802697120^3 & = \\
1117386592077753452^3 & + & 14123302420417013824^3 & = \\
1497369344185092651^3 & + & 14120022667932733461^3 & = \\
2223357078845220136^3 & + & 14107248762203982476^3 & = \\
4163116835733008647^3 & + & 14004053464077523769^3 & = \\
6716379921779399326^3 & + & 13600192974314732786^3 & = \\
7925282888762885516^3 & + & 13239637283805550696^3& = \\
8548057588027946352^3 & + & 12993955521710159724^3& = \\
8937735731741157614^3 & + & 12815060285137243042^3& = \\
11005214445987377356^3 & + & 11410505395325664056^3&  \\
\end{array}
\]

\noindent$T_{12}=2971^3\,T_{11}=\\
73914858746493893996583617733225161086864012865017882136931801625152$:
\[
\begin{array}{rcrc}
845205202844653597674^3  & + & 41967142660804626363462^3 & = \\
1933097542618122241026^3 & + & 41965889731136229476526^3 & = \\
1952714722754103222628^3 & + & 41965847682542813143520^3 & = \\
3319755565063005505892^3 & + & 41960331491058948071104^3 & = \\
4448684321573910266121^3 & + & 41950587346428151112631^3 & = \\
6605593881249149024056^3 & + & 41912636072508031936196^3  & = \\
12368620118962768690237^3 & + & 41606042841774323117699^3 & = \\
19954364747606595397546^3 & + & 40406173326689071107206^3 & = \\
23546015462514532868036^3 & + & 39334962370186291117816^3 & = \\
25396279094031028611792^3 & + & 38605041855000884540004^3 & = \\
26554012859002979271194^3 & + & 38073544107142749077782^3 & = \\
32696492119028498124676^3 & + & 33900611529512547910376^3 &  \\
\end{array}
\]

\noindent$T_{13}=4327^3\,T_{12}=\\
5988146776742829080553965820313279739849705084894534523771076163371248442670016$:
\[
\begin{array}{rcrc}

3657202912708816117135398^3 & + & 181591826293301618274700074^3  & = \\
8364513066908614936919502^3 & + & 181586404866626464944928002^3  & = \\
8449396605357004644311356^3 & + & 181586222922362752472011040^3  & = \\
14364582330027624823994684^3 & + & 181562354361812068303667008^3 & = \\
19249457059450309721505567^3 & + & 181520191447994609864354337^3 & = \\
28582404724165067827090312^3 & + & 181355976285742254187920092^3 & = \\
53519019254751900122655499^3 & + & 180029347376357496130283573^3 & = \\
54818831102057750995052604^3 & + & 179911586979069103444414128^3 & = \\
86342536262893738285181542^3 & + & 174837511984583610680880362^3 & = \\
101883608906300383719991772^3 & + & 170202382175796081666789832^3& = \\
109889699639872260803223984^3 & + & 167044016106588827404597308^3& = \\
114899213640905891306456438^3 & + & 164744225351606675259562714^3& = \\
141477721399036311385473052^3 & + & 146687946088200794808196952^3&  \\
\end{array}
\]

\noindent$T_{14}=7549^3\,T_{13}=$
\[
\begin{array}{c}
257608810925730001281963766003343299028977072     ~~~\backslash\\
~~~~~~~~~~~~~~~~~~~
5881505682307757452553496715044742867424072384:
\end{array}
\]
\[
\begin{array}{rcrc}
27608224788038852868255119502^3 & + & 1370836696688133916355710858626^3 & = \\
63143709142093134158805320598^3 & + & 1370795770338163183869261487098^3 & = \\
63784494973840028059906426444^3 & + & 1370794396840916418411211340960^3 & = \\
108438232009378539796335869516^3 & + & 1370614213077319303624382243392^3 & = \\
145314151341790388087645525283^3 & + & 1370295925240911309866010890013^3 & = \\
215768573262722097026704765288^3 & + & 1369056264981068276864608774508^3 & = \\
404015076354122094025926361951^3 & + & 1359041543344122738287510692577^3 & = \\
413827355989433962261652107596^3 & + & 1358152570104992661901882252272^3 & = \\
617989830682279948575932296880^3 & + & 1327627770274178602420131034444^3 & = \\
651799806248584830314835460558^3 & + & 1319848377971621677029965852738^3 & = \\
769119363633661596702217886828^3 & + & 1284857783045084620502596441768^3 & = \\
829557342581395696803537855216^3 & + & 1261015277588639058077305078092^3 & = \\
867374163775198573472439650462^3 & + & 1243654157179278791534438927986^3 & = \\
1068015318841325114648936069548^3 & + & 1107347305019827800007078790648^3 &  \\
\end{array}
\]


\end{document}